\magnification=1200

\def\Q{{\bf {Q}}}

\def\Z{{\bf Z}}

  \def\cN{{\cal {N}}}

%
\catcode`@=11
%
%
\def\bibn@me{R\'ef\'erences}
\def\bibliographym@rk{\centerline{{\sc\bibn@me}}
	\sectionmark\section{\ignorespaces}{\unskip\bibn@me}
	\bigbreak\bgroup
	\ifx\ninepoint\undefined\relax\else\ninepoint\fi}
%
%
%
\let\refsp@ce=\ 
\let\bibleftm@rk=[
\let\bibrightm@rk=]
%
%
%
\def\numero{n\raise.82ex\hbox{$\fam0\scriptscriptstyle o$}~\ignorespaces}
%
%
\newcount\equationc@unt
\newcount\bibc@unt
\newif\ifref@changes\ref@changesfalse
\newif\ifpageref@changes\ref@changesfalse
\newif\ifbib@changes\bib@changesfalse
\newif\ifref@undefined\ref@undefinedfalse
\newif\ifpageref@undefined\ref@undefinedfalse
\newif\ifbib@undefined\bib@undefinedfalse
\newwrite\@auxout
%
%
\def\eqnum{\global\advance\equationc@unt by 1%
\edef\lastref{\number\equationc@unt}%
\eqno{(\lastref)}}
%
%
%
%
%
%
\def\re@dreferences#1#2{{%
	\re@dreferenceslist{#1}#2,\undefined\@@}}
\def\re@dreferenceslist#1#2,#3\@@{\def\next{#2}%
	\expandafter\ifx\csname#1@@\meaning\next\endcsname\relax
	??\immediate\write16
	{Warning, #1-reference "\next" on page \the\pageno\space
	is undefined.}%
	\global\csname#1@undefinedtrue\endcsname
	\else\csname#1@@\meaning\next\endcsname\fi
	\ifx#3\undefined\relax
	\else,\refsp@ce\re@dreferenceslist{#1}#3\@@\fi}
%
%
%
\def\newlabel#1#2{{\def\next{#1}\newl@bel#2}}
\def\newl@bel#1#2{%
	\expandafter\xdef\csname ref@@\meaning\next\endcsname{#1}%
	\expandafter\xdef\csname pageref@@\meaning\next\endcsname{#2}}
\def\label#1{{%
	\toks0={#1}\message{ref(\lastref) \the\toks0,}%
	\ignorespaces\immediate\write\@auxout%
	{\noexpand\newlabel{\the\toks0}{{\lastref}{\the\pageno}}}%
	\def\next{#1}%
	\expandafter\ifx\csname ref@@\meaning\next\endcsname\lastref%
	\else\global\ref@changestrue\fi%
	\newlabel{#1}{{\lastref}{\the\pageno}}}}
\def\ref#1{\re@dreferences{ref}{#1}}
\def\pageref#1{\re@dreferences{pageref}{#1}}
%
%
\def\bibcite#1#2{{\def\next{#1}%
	\expandafter\xdef\csname bib@@\meaning\next\endcsname{#2}}}
\def\cite#1{\bibleftm@rk\re@dreferences{bib}{#1}\bibrightm@rk}
%
%
\def\beginthebibliography#1{\bibliographym@rk
	\setbox0\hbox{\bibleftm@rk#1\bibrightm@rk\enspace}
	\parindent=\wd0
	\global\bibc@unt=0
	\def\bibitem##1{\global\advance\bibc@unt by 1
		\edef\lastref{\number\bibc@unt}
		{\toks0={##1}
		\message{bib[\lastref] \the\toks0,}%
		\immediate\write\@auxout
		{\noexpand\bibcite{\the\toks0}{\lastref}}}
		\def\next{##1}%
		\expandafter\ifx
		\csname bib@@\meaning\next\endcsname\lastref
		\else\global\bib@changestrue\fi%
		\bibcite{##1}{\lastref}
		\medbreak
		\item{\hfill\bibleftm@rk\lastref\bibrightm@rk}%
		}
	}
\def\endthebibliography{\egroup\par}
%
%
%
\def\@closeaux{\closeout\@auxout
	\ifref@changes\immediate\write16
	{Warning, changes in references.}\fi
	\ifpageref@changes\immediate\write16
	{Warning, changes in page references.}\fi
	\ifbib@changes\immediate\write16
	{Warning, changes in bibliography.}\fi
	\ifref@undefined\immediate\write16
	{Warning, references undefined.}\fi
	\ifpageref@undefined\immediate\write16
	{Warning, page references undefined.}\fi
	\ifbib@undefined\immediate\write16
	{Warning, citations undefined.}\fi}
%
%
\immediate\openin\@auxout=\jobname.aux
\ifeof\@auxout \immediate\write16
  {Creating file \jobname.aux}
\immediate\closein\@auxout
\immediate\openout\@auxout=\jobname.aux
\immediate\write\@auxout {\relax}%
\immediate\closeout\@auxout
\else\immediate\closein\@auxout\fi
%
%
\input\jobname.aux
\immediate\openout\@auxout=\jobname.aux
%
%
\catcode`@=12

%
\catcode`@=11
\def\bibliographym@rk{\bgroup}
%
%
\outer\def\bye{ 	\par\vfill\supereject\end}

\def\house#1{\setbox1=\hbox{$\,#1\,$}%
\dimen1=\ht1 \advance\dimen1 by 2pt \dimen2=\dp1 \advance\dimen2 by 2pt
\setbox1=\hbox{\vrule height\dimen1 depth\dimen2\box1\vrule}%
\setbox1=\vbox{\hrule\box1}%
\advance\dimen1 by .4pt \ht1=\dimen1
\advance\dimen2 by .4pt \dp1=\dimen2 \box1\relax}

  \def\eps{{\varepsilon}}

  \def\noi{\noindent}

\def\build#1_#2^#3{\mathrel{\mathop{\kern 0pt#1}\limits_{#2}^{#3}}}

\def\date {le\ {\the\day}\ \ifcase\month\or janvier
\or fevrier\or mars\or avril\or mai\or juin\or juillet\or
ao\^ut\or septembre\or octobre\or novembre
\or d\'ecembre\fi\ {\oldstyle\the\year}}

\font\fivegoth=eufm5 \font\sevengoth=eufm7 \font\tengoth=eufm10

\newfam\gothfam \scriptscriptfont\gothfam=\fivegoth
\textfont\gothfam=\tengoth \scriptfont\gothfam=\sevengoth

\def\smallsquare{\vbox{\hrule\hbox{\vrule height 1 ex\kern 1 ex\vrule}\hrule}}
\def\cqfd{\hfill \smallsquare\vskip 3mm}

\def\og{\leavevmode\raise.3ex\hbox{$\scriptscriptstyle 
\langle\!\langle\,$}}
\def \fg {\leavevmode\raise.3ex\hbox{$\scriptscriptstyle 
\!\rangle\!\rangle\,\,$}}

\def\rme{{\rm e}}

\def\bfv{{\bf v}}
\def\bfb{{\bf b}}
\def\bfx{{\bf x}}
\def\bfX{{\bf X}}
\def\bfn{{\bf n}}

\def\cN{{\cal N}}

\centerline{}

\vskip 8mm

\centerline{\bf On the digital representation of integers with bounded prime factors}

\vskip 13mm

\centerline{Yann B{\sevenrm UGEAUD} \footnote{}{\rm
}}

{\narrower\narrower
\vskip 15mm

\proclaim Abstract. {
Let $b \ge 2$ be an integer. 
Not much is known on the representation in base $b$ of prime numbers or of 
numbers whose prime factors belong to a given, finite set.  
Among other results, we establish that any sufficiently large integer which is not a multiple of $b$
and has only small (in a suitable sense) prime factors has at least four nonzero digits in 
its representation in base $b$. 
}

}

\vskip 15mm

\centerline{\bf 1. Introduction and results}

\vskip 5mm

We still do not know whether there are infinitely many prime numbers of the form
$2^n + 1$ (that is, with only two nonzero binary digits) or of the form
$11 \ldots 11$ (that is, with only the digit $1$ in their decimal representation). 
Both questions are notorious, very difficult open problems, which at present seem to be 
completely out of reach. 
However, there have been recently several spectacular advances on the digital 
representation of prime numbers. 
In 2010, Mauduit and Rivat \cite{MauRi10} established that 
the sum of digits of primes is well-distributed.
Subsequently, Bourgain \cite{Bour15} showed the existence of prime
numbers in the sparse set defined by prescribing a positive proportion of the binary digits.
This year, Maynard \cite{May16} proved that, if $d$ is any digit in $\{0, 1, \ldots , 9\}$, then 
there exist infinitely many prime numbers which do not have the digit $d$ in 
their decimal representation. The proofs of all these results depend largely on
Fourier analysis techniques. 

In the present note we study a related problem, namely the digital representation of 
integers all of whose prime factors belong to a finite, given 
set $S$ of prime numbers.
We apply techniques from Diophantine approximation 
to discuss the following general (and left intentionally vague) question: 

\smallskip
{\it Do there exist arbitrarily large integers 
which have only small prime factors and, at the same time, few nonzero digits in their
representation in some integer base?}

\smallskip

The expected answer is {\it no} and our results are a modest step in this direction.

Let $n$ be a positive integer $n$ and $P[n]$ denote its greatest prime factor, with the 
convention that $P[1] = 1$. 
Let $S = \{q_1, \ldots , q_s\}$ be a finite, non-empty set of distinct prime numbers.
Write $n = q_1^{r_1} \ldots q_s^{r_s} M$, where 
$r_1, \ldots , r_s$ are non-negative integers and $M$ is an integer 
relatively prime to $q_1 \ldots q_s$. We define the $S$-part $[n]_S$ 
of $n$ by 
$$
[n]_S := q_1^{r_1} \ldots q_s^{r_s}.
$$
The $S$-parts of linear 
recurrence sequences and of integer polynomials and decomposable forms 
evaluated at integer points have been studied in \cite{GrVi13,BuEv16,BuEvGy16}. 

In the sequel, for a given integer $k \ge 2$, we denote by $(u_j^{(k)})_{j \ge 1}$ 
the sequence, arranged in increasing order, of all positive integers which are 
not divisible by $b$ and have at most $k$ nonzero digits in their $b$-ary representation. 
Said differently, $(u_j^{(k)})_{j \ge 1}$ is the ordered sequence composed of the integers of the
form
$$
d_k b^{n_k} + \ldots + d_2 b^{n_2} + d_1, 
\quad n_k > \ldots > n_2 > 0, \quad
d_1, \ldots , d_k \in \{0, 1, \ldots , b-1\}, \quad
d_1 d_k \not= 0.
$$
We stress that, for the questions investigated in the present note, 
it is natural to restrict 
our attention to integers not divisible by $b$. 
Obviously, the sequence $(u_j^{(k)})_{j \ge 1}$ depends on $b$, but, for 
shortening the notation, we have decided not to mention this dependence.

Our first result shows that, for any base $b$, 
there are only finitely many integers not divisible by $b$ which have a given number of 
nonzero $b$-ary digits and whose prime divisors belong to a given finite set.

\proclaim Theorem 1.1. 
Let $b \ge 2, k \ge 2$ be integers and $\eps$ a positive real number.
Let $S$ be a finite, non-empty set of prime numbers. 
Then, we have 
$$
[u_j^{(k)}]_S <  (u_j^{(k)})^{\eps}, 
$$
for every sufficiently large integer $j$. 
In particular, the greatest prime factor of $u_j^{(k)}$ tends to infinity as $j$ tends to
infinity.

The proof of Theorem 1.1 rests on the Schmidt Subspace Theorem and does not allow us
to estimate the speed with which $P[u_j^{(k)}]$ tends to infinity with $j$. 
It turns out that, by means of the theory of linear forms in logarithms, we are able to derive
such an estimate, but (apparently) only for $k \le 3$. 

The case $k=2$ has already been considered. It reduces to the study of a finite 
union of binary linear recurrences of the form 
$$
(d_2 b^n + d_1 1^n)_{n \ge 1}, \quad \hbox{where 
$d_1, d_2$ are digits in $\{1, \ldots , b-1\}$}.
$$ 
We gather in the next theorem a recent result of Bugeaud and Evertse \cite{BuEv16} 
and an immediate consequence of a lower bound for the greatest prime 
factor of terms of binary recurrence sequences, established by Stewart \cite{Ste13b}. 

\proclaim Theorem BES.
Let $b \ge 2$ be an integer. 
Let $S$ be a finite, non-empty set of prime numbers. 
Then, there exist an effectively computable positive number $c_1$, depending 
only on $b$, and an effectively computable positive number $c_2$, depending 
only on $b$ and $S$, such that
$$
[u_j^{(2)}]_S \le (u_j^{(2)})^{1 - c_1}, \quad
\hbox{for every $j \ge c_2$}.
$$
Furthermore, there exists an effectively computable positive number $c_3$, depending 
only on $b$ and $S$, such that
$$
P[u_j^{(2)}] 
> (\log u_j^{(2)})^{1/2} \,  \exp \Bigl( {\log \log u_j^{(2)} \over 105 \log \log \log u_j^{(2)}} \Bigr),
 \quad \hbox{for $j > c_3$}.
$$

We point out that the constant $c_1$ in Theorem BES does not depend on $S$.

The main new result of the present note is an 
estimate of the speed with which $P[u_j^{(3)}]$ tends to infinity with $j$.

\proclaim Theorem 1.2. 
Let $b \ge 2$ be an integer. 
Let $S$ be a finite, non-empty set of prime numbers. 
Then, there exist effectively computable positive numbers $c_4$ and $c_5$, depending 
only on $b$ and $S$, such that
$$
[u_j^{(3)}]_S \le (u_j^{(3)})^{1 - c_4}, \quad
\hbox{for every $j \ge c_5$}.
$$
Furthermore, for every positive real number $\eps$, there exists an effectively
computable positive number $c_6$, depending 
only on $b$ and $\eps$, such that 
$$
P[u_j^{(3)}] > (1 - \eps) \log \log u_j^{(3)} 
\, {\log \log \log u_j^{(3)} \over \log \log \log \log u_j^{(3)}}, \quad
\hbox{for $j > c_6$}.   \eqno (1.1) 
$$

The proof of Theorem 1.2 yields a very small admissible value for $c_4$.

We point out the following reformulation of the second assertion of Theorem 1.2.
Recall that a positive integer is called $B$-smooth if all its prime factors are less than 
or equal to $B$. 

\proclaim Corollary 1.3.
Let $b \ge 2$ be an integer.
Let $\eps$ be a positive integer.
There exists an effectively computable positive integer $n_0$, depending 
only on $b$ and $\eps$, such that any integer $n > n_0$ which is 
not divisible by $b$ and is 
$$
(1 - \eps) (\log \log n) {\log \log \log n \over \log \log \log \log n}\hbox{-smooth}
$$
has at least four nonzero digits in its $b$-ary representation. 

It is very likely that any large integer cannot be `very' smooth  
and, simultaneously, have only few nonzero digits in its $b$-ary representation. 
Corollary 1.3 provides a first result in this direction. 

The proofs of our theorems are obtained by direct applications of classical deep tools 
of Diophantine approximation, namely the Schmidt Subspace Theorem and the theory 
of linear forms in the logarithms of algebraic numbers. The latter theory 
has already been applied to get lower bounds for the greatest 
prime factor of linear recurrence sequences (under some assumptions, see \cite{Ste08c}) and 
for the greatest prime factor of integer polynomials and
decomposable forms evaluated at integer points 
(see e.g. \cite{GyYu06}). The bounds obtained in \cite{Ste08c,GyYu06} 
have exactly the same 
order of magnitude as our bound in Theorem 1.2, 
that is, they involve a double logarithm times a triple logarithm
divided by a quadruple logarithm. A brief explanation is given at the end of Section 2. 

An interesting feature of the proof of Theorem 1.2 is that it combines 
estimates for Archimedean and non-Archimedean linear forms in logarithms. 
Similar arguments appeared when searching for perfect powers with 
few digits; see \cite{BeBuMi12,BeBuMi13}.


\vskip 5mm

\centerline{\bf 2. Auxiliary results from Diophantine approximation}

\vskip 5mm

The Schmidt Subspace Theorem \cite{Schm70a,Schm72,SchmLN} 
is a powerful multidimensional extension of the Roth Theorem. 
We quote below a version of it
which is suitable for our purpose, but the reader should
keep in mind that there are more general
formulations.

\proclaim Theorem 2.1. 
Let $m\ge 2$ be an integer.
Let $S'$ be a finite set of prime numbers.
Let $L_{1, \infty}, \ldots , L_{m, \infty}$ 
be $m$ linearly independent linear forms in $m$ variables with integer coefficients.
For any prime $\ell$ in $S'$, let $L_{1, \ell}, \ldots , L_{m, \ell}$ 
be $m$ linearly independent linear forms in $m$ variables with integer coefficients.
Let $\eps$ be a positive real number.
Then, there are an integer $T$ and proper subspaces $S_1, \ldots , S_T$
of $\Q^m$ such that all the solutions $\bfx = (x_1, \ldots, x_m)$
in $\Z^m$ to the inequality
$$ 
\prod_{\ell \in S'} \, \prod_{i=1}^m \,
\vert L_{i, \ell} (\bfx) \vert_{\ell}
\cdot \prod_{i=1}^m \,
{\vert L_{i, \infty} (\bfx) \vert} \,\le \, 
(\max\{1, |x_1|, \ldots , |x_m|\})^{-\eps}  
$$
are contained in the union $S_1 \cup \ldots \cup S_T$.

We quote an immediate 
corollary of a theorem of Matveev \cite{Matv00}.

\proclaim Theorem 2.2.  
Let $n \ge 2$ be an integer. 
Let $x_1/y_1, \ldots, x_n/y_n$ be positive rational numbers. 
Let $b_1, \ldots, b_n$ be integers such that $(x_1/y_1)^{b_1} \ldots (x_n/y_n)^{b_n} \not= 1$. 
Let $A_1, \ldots, A_n$ be real numbers with
$$
A_i \ge \max\{|x_i|, |y_i|, \rme\}, \quad 1\le i \le n.
$$
Set
$$
B = \max\Bigl\{1, \max\Bigl\{ |b_j| \ {\log A_j \over \log A_n} : 1 \le j \le n \Bigr\} \Bigr\}.
$$
Then, we have
$$
\log \Bigl| \Bigl( {x_1 \over y_1} \Bigr)^{b_1} \ldots \Bigl( {x_n \over y_n} \Bigr)^{b_n}  - 1 \Bigr|  
> - 8 \times 30^{n+3} \, n^{9/2}  \,  \log (\rme B) \, 
\log A_1 \ldots \log A_n.       
$$

The next statement was proved by Yu \cite{Yu07}. 
For a prime number $p$ and a nonzero rational number $z$ we denote by $v_p(z)$ the exponent 
of $p$ in the decomposition of $z$ in product of prime factors. 

\proclaim Theorem 2.3.
Let $n \ge 2$ be an integer. 
Let $x_1/y_1, \ldots, x_n/y_n$ be nonzero rational numbers. 
Let $b_1, \ldots, b_n$ be nonzero 
integers such that $(x_1/y_1)^{b_1} \ldots (x_n/y_n)^{b_n} \not= 1$. 
Let $B$ and $B_n$ be real numbers such that
$$
B \ge \max\{|b_1|, \ldots , |b_n|, 3\} \quad 
\hbox{and} \quad
B \ge B_n \ge |b_n|.
$$
Assume that
$$
v_p (b_n) \le v_p (b_j), \quad j = 1, \ldots , n.
$$
Let $\delta$ be a real number with $0< \delta \le 1/2$. 
Then, we have
$$
\eqalign{
v_p \Bigl( \Bigl( {x_1 \over y_1} \Bigr)^{b_1} \ldots \Bigl( {x_n \over y_n} \Bigr)^{b_n}  - 1 \Bigr)  < & \, 
(16 \rme)^{2(n+1)} n^{3/2}  \, (\log (2n))^2
 \,  {p \over (\log p)^2}  \cr
&  \, \, \, \, 
\max\Bigl\{  (\log A_1) \cdots (\log A_n) (\log T), {\delta B \over B_n} \Bigr\}, \cr}
$$
where
$$
T = 2 B_n \delta^{-1}  \rme^{(n+1)(6n+5)} p^{n+1}  (\log A_1) \cdots (\log A_{n-1}). 
$$

There are two key ingredients in Theorems 2.2 and 2.3 
which explain the quality of the estimates in Theorem 1.2.
A first one is the dependence on $n$, which is only exponential: this allows us to get 
the extra factor triple logarithm over quadruple logarithm. The use of earlier estimates for
linear forms in logarithms would give only the factor involving the double logarithm in (1.1). 
A second one is the factor $\log A_n$ occurring 
in the denominator in the
definition of $B$ in the statement of Theorem 2.2. 
The formulation of Theorem 2.3 is slightly different, but, in our special case, 
it yields a similar refinement.  This allows us to save a (small) power of $u_j^{(3)}$ 
when estimating its $S$-part. Without this refinement, the saving would 
be much smaller, namely less than any power of $u_j^{(3)}$.

\vskip 5mm

\centerline{\bf 3. Proofs}

\vskip 5mm

\noindent {\it Proof of Theorem 1.1.}  

Let $k \ge 2$ be an integer and $\eps$ a positive real number.
Let $\cN_1$ be the set of $k$-tuples 
$(n_k, \ldots , n_2, n_1)$ such that 
$n_k > \ldots > n_2 > n_1 = 0$ and 
$$
[d_k b^{n_k} + \cdots + d_2 b^{n_2} + d_1]_S > 
(d_k b^{n_k} + \cdots + d_2 b^{n_2} + d_1)^{\eps},
$$ 
for some integers $d_1, \ldots , d_k$ in $\{0, \ldots , b-1\}$ 
such that $d_1 d_k \not= 0$.

Assume that $\cN_1$ is infinite. 
Then, there exist an integer $h \ge 2$, positive integers $D_1, \ldots , D_h$,
an infinite set $\cN_2$, contained in $\cN_1$, of $h$-tuples 
$(n_{h,i}, \ldots , n_{1,i})$ such that
$n_{h,i} > \ldots > n_{1,i} \ge 0$,
$$
[D_h b^{n_{h,i}} + \cdots + D_1 b^{n_{1,i}}]_S > 
(D_h b^{n_{h,i}} + \cdots + D_1 b^{n_{1,i}})^{\eps}, \quad i \ge 1, 
$$
and
$$
\lim_{i \to + \infty} \, (n_{\ell,i} - n_{\ell-1,i}) = + \infty, \quad \ell = 2, \ldots , h.   \eqno (3.1)
$$
We are in position to apply Theorem 2.1.

Let $S_1$ denote the set of prime divisors of $b$. Without any loss of generality, we may assume
that $S$ and $S_1$ are disjoint. 
Consider the linear forms in $\bfX = (X_1, \ldots , X_h)$ given by 
$$
L_{j, \infty} (\bfX) := X_j, \quad  1 \le j \le h, \quad
$$
and, for every prime number $p$ in $S_1$,
$$
 L_{j, p} (\bfX) := X_j, \quad 1 \le j \le h, \quad
$$
and, for every prime number $p$ in $S$,
$$
 L_{j, p} (\bfX) := X_j, \quad 1 \le j \le h-1, \quad
L_{h, p} (\bfX) := D_h X_h + \ldots + D_1 X_1.
$$
By Theorem 2.1 applied with $S' = S \cup S_1$,
the set of tuples $\bfb = (b^{n_h}, \ldots , b^{n_2}, b^{n_1})$ 
such that $n_h > \ldots > n_1 \ge 0$ and 
$$
\prod_{j=1}^h \, |L_{j, \infty} (\bfb)| 
\times \prod_{p \in S \cup S_1} \,
\prod_{j=1}^h \, |L_{j, p} (\bfb)|_p < b^{-\eps n_h}    \eqno (3.2)
$$
is contained in a finite union of proper subspaces of $\Z^h$. 
Since the left hand side of (3.2) is equal to 
$[D_h b^{n_h} + \cdots + D_1 b^{n_1}]_S^{-1}$, the set of tuples $(b^{n_h}, \ldots , b^{n_1})$,
where $(n_h, \ldots , n_1)$ lies in $\cN_2$, is contained in a 
finite union of proper subspaces of $\Z^h$. 

Thus, there exist integers $t_1, \ldots , t_h$, not all zero, and an infinite
set $\cN_3$, contained in $\cN_2$, 
of integer tuples $(n_h, \ldots , n_1)$ such that $n_h > \ldots > n_1 \ge 0$ and 
$$
t_h b^{n_h} + \cdots + t_1 b^{n_1} = 0.
$$
We then deduce from (3.1) that $t_1 = \ldots = t_h = 0$, a contradiction. 
Consequently, the set $\cN_1$ must be finite. This establishes the theorem. 
\cqfd

\vskip 5mm

\noindent {\it Proof of Theorem 1.2.}  

Below, the constants $c_1, c_2, \ldots$ are effectively computable and depend at most 
on $b$ and the constants $C_1, C_2, \ldots$ are absolute and effectively computable. 

Let $q_1, \ldots , q_s$ be distinct prime numbers written in increasing order. 
Let $j \ge b^4$ be an integer and write
$$
u_j^{(3)} = d_3 b^m + d_2 b^n + d_1, 
\quad \hbox{where $d_1, d_2, d_3  \in \{0, 1, \ldots , b-1\}$,
$d_1 d_3 \not= 0$, $m  > n > 0$.}
$$
There exist non-negative integers $r_1, \ldots , r_s$ and a positive integer $M$ 
coprime with $q_1 \ldots q_s$ such that
$$
u_j^{(3)}  = q_1^{r_1} \cdots q_s^{r_s} M.
$$
Assume first that $m \ge 2 n$. 
Since
$$
\Lambda_a := |q_1^{r_1} \cdots q_s^{r_s} b^{-m} (M d_3^{-1})  - 1| 
\le  b^{1+n-m} \le  b^{-(m-2)/2},
$$
we get the upper bound 
$$
\log \Lambda_a  \le -  \Bigl({m \over 2} - 1 \Bigr) \, \log b. 
$$
For the lower bound, setting
$$
Q := (\log q_1) \cdots (\log q_s) \quad 
\hbox{and} \quad A := \max\{M, d_3, 2\}, 
$$
and using that $r_j \log q_j \le (m+1) \log b$ for $j = 1, \ldots , s$, 
Theorem 2.2 implies that
$$
\log \Lambda_a \ge - c_1 C_1^s \, Q \, (\log A) \, \log {m  \over \log A}.
$$
Comparing both estimates, we deduce that
$$
m  \le c_2 \,  C_2^s \, Q \,  (\log Q) \, (\log A).   \eqno (3.3)
$$

Assume now that $m \le 2n$.
Let $p$ be the smallest prime divisor of $b$. Set 
$$
\Lambda_u := q_1^{r_1} \cdots q_s^{r_s} {M \over d_1}  - 1 
= {b^n \over d_1} (d_2 + d_3 b^{m-n})
$$
and
$$
A = \max\{M, d_1, 2\}, \quad B = \max\{r_1, \ldots , r_s, 3\}. 
$$
Observe that 
$$
v_p (\Lambda_u) \ge n - {\log b \over \log p} \ge {m \over 2} - {\log b \over \log p} .      \eqno (3.4)
$$
It follows from Theorem 2.3 applied with 
$$
\delta = {Q (\log A) \over B}
$$
that
$$
B < 2 Q \log A, \quad \hbox{if $\delta > 1/2$},   \eqno (3.5)
$$
and, otherwise,
$$
v_p (\Lambda_u) < c_3 C_3^s Q (\log A) \max\Bigl\{ \log \Bigl( {B \over \log A} \Bigr), 1 \Bigr\} 
\eqno (3.6)
$$
Observe that
$$
c_4 m \le B \le c_5 m.   \eqno (3.7)
$$
We deduce from (3.5) and (3.7) that
$$
m < c_6 Q \log A,    \eqno (3.8) 
$$
and, by combining (3.4), (3.6), and (3.7), we get
$$
m  \le c_7 \,  C_4^s \, Q \,  (\log Q) \, (\log A).  \eqno (3.9)
$$

Observe that $A \le \max\{M, b\}$. 
It then follows from (3.3), (3.8), and (3.9) that if
$$
m > c_8 \,  C_5^s \, Q \, (\log Q) \, (\log b),
$$
then $A = M$ and, using that $m \log b \le \log u_j^{(3)}$, we conclude that 
$$
M \ge (u_j^{(3)})^{(c_9 C_5^s \, Q \,  (\log Q))^{-1}},
$$
thus
$$
[u_j^{(3)}]_S = {u_j^{(3)} \over M} \le (u_j^{(3)})^{1 - (c_9 C_5^s \, Q \,  (\log Q))^{-1}}.
$$

In the particular case where $M=1$ and $q_1, \ldots , q_s$ 
are the first $s$ prime numbers $p_1, \ldots , p_s$, written
in increasing order, the above proof shows that
$$
m \le c_{10} \, C_6^s \, \bigl( \prod_{k=1}^s \log p_k \bigr) \, (\log \log p_s).  \eqno (3.10)
$$
Let $\eps$ be a positive real number. 
We deduce from (3.10), the Prime Number Theorem, and the inequality
$\log u_j^{(3)} < (m+1) \log b$ that
$$
\log \log u_j^{(3)} \le c_{11} + C_7 s + \sum_{k=1}^s \log \log p_k  
\le (1+ \eps) p_s {\log \log p_s \over \log p_s},
$$
if $j$ is sufficiently large in terms of $\eps$. 
This establishes (1.1) and completes the proof of Theorem 1.2. 
\cqfd

\vskip 5mm

\centerline{\bf 4. Additional remarks}

\vskip 5mm

In this section, we present additional results, discuss related problems, 
and make some suggestions for further research.

Arguing as Stewart did in \cite{Ste08c}, 
we can apply the arithmetic-geometric mean inequality in the course of the
proof of Theorem 1.2 
to derive a lower bound for $Q[u_j^{(3)}]$, where $Q[n]$ denotes the
greatest square-free divisor of a positive integer $n$.

\proclaim Theorem 4.1. 
Let $b \ge 2$ be an integer. 
There exist effectively
computable positive numbers $c_1, c_2$, depending 
only on $b$, such that 
$$
Q[u_j^{(3)}] >  \exp \, \Bigl( c_1 \,  \log \log u_j^{(3)} 
\, {\log \log \log u_j^{(3)} \over \log \log \log \log u_j^{(3)}} \, \Bigr), \quad
\hbox{for $j > c_2$}. 
$$

It is not difficult to make Theorem 1.2 completely explicit. Even in the special case where 
the cardinality of the set $S$ is small, the bounds obtained are rather large, since 
estimates for linear forms in three or more logarithms are needed. 
Thus, it is presumably not straightforward to solve completely an equation like
$$
2^a + 2^b + 1 = 3^x 5^y, \quad
\hbox{in non-negative integers $a, b, x, y$ with $a > b > 0$.}
$$

Presumably, other techniques, based on the hypergeometric method, 
could yield effective improvements of Theorem 1.2 in some special cases. 

\proclaim Problem 4.2.
Take $S = \{3, 5\}$. 
Prove that there exists an effectively computable integer $m_0$ such that,
for any integers $m, n$ with $m > m_0$ and $m > n > 0$, we have
$$
[2^m + 2^n + 1]_S \le 2^{3m/4}.     \eqno (4.1)
$$

No importance should be attached to the value $3/4$ in (4.1). 
Similar questions have been successfully addressed in \cite{BeFiTr08,BeFiTr09}.

Perfect powers with few nonzero digits in some given integer base have been 
studied in \cite{CoZa00,CoZa12,Ben12,BeBuMi12,BeBuMi13,BeBu14}; 
see also the references given therein. 
We briefly discuss a related problem.  
Let $a, b$ be integers 
such that $a > b > 1$. Perfect powers in the bi-infinite
sequence $(a^m + b^n + 1)_{m, n \ge 1}$ have been considered by 
Corvaja and Zannier \cite{CoZa05} and also in \cite{BeBuMi12,BeBuMi13}. 
All the general results obtained so far have been established 
under the assumption that 
$a$ and $b$ are not coprime. To remove 
this coprimeness assumption seems to be a very difficult problem.

We note that the methods of the proof of Theorem 1.2 allows us to establish the 
following result.

\proclaim Theorem 4.3.
Let $a, b$ be distinct integers with $\gcd (a, b) \ge 2$. 
Let $\bfv = (v_j)_{j \ge 1}$ denote the increasing sequence composed of all the 
integers of the form $a^m + b^n + 1$, with $m, n \ge 1$. Then,
for every positive $\eps$, we have
$$
P [v_j] > (1 - \eps) \log \log v_j \, {\log \log \log v_j \over \log \log \log \log v_j},
$$
when $j$ exceeds some effectively computable constant depending only on $a$ and $b$.

We point out the following problem, which is probably rather difficult.

\proclaim Problem 4.4.
Give an effective lower bound for the greatest prime factor of 
$2^m + 3^n + 1$ in terms of $\max\{m, n\}$.

\vskip 7mm

\centerline{\bf References}

\vskip 7mm

\beginthebibliography{999}

\medskip 

\bibitem{Ben12}
M. A. Bennett, 
{\it Perfect powers with few ternary digits},
Integers 12A (Selfridge Memorial Volume), 8pp., 2012. 

\bibitem{BeBu14}
M. A. Bennett and Y. Bugeaud,
{\it Perfect powers with three digits}, 
Mathematika 60 (2014),  66--84.

\bibitem{BeBuMi12}
M. A. Bennett, Y. Bugeaud and M. Mignotte,
{\it Perfect powers with few binary digits
and related Diophantine problems, II},
Math. Proc. Cambridge Philos. Soc. 153 (2012), 525--540.

\bibitem{BeBuMi13}
M. A. Bennett, Y. Bugeaud and M. Mignotte,
{\it Perfect powers with few binary digits
and related Diophantine problems},
 Ann. Sc. Norm. Super. Pisa Cl. Sci. 12 (2013), 525--540.

\bibitem{BeFiTr08}
M. A. Bennett, M. Filaseta, and O. Trifonov, 
{\it Yet another generalization of the Ramanujan-Nagell equation}, 
Acta Arith. 134 (2008), 211--217.

\bibitem{BeFiTr09}
M. A. Bennett, M. Filaseta, and O. Trifonov, 
{\it On the factorization of consecutive integers}, 
J. Reine Angew. Math. 629 (2009), 171--200.

\bibitem{Bour15}
J. Bourgain,
{\it Prescribing the binary digits of primes, II},
Israel J. Math. 206 (2015), 165--182.

\bibitem{BuEv16}
Y. Bugeaud and J.-H. Evertse, 
{\it On the $S$-part of integer recurrent sequences}.
In preparation. 

\bibitem{BuEvGy16}
Y. Bugeaud, J.-H. Evertse, and K. Gy\H ory, 
{\it On $S$-parts of decomposable forms at integral points}.
In preparation. 

\bibitem{CoZa00}
P. Corvaja and U. Zannier,
{\it On the Diophantine equation $f(a^m,y)=b^n$},
Acta Arith.  94  (2000),  25--40.

\bibitem{CoZa05}
P. Corvaja and U. Zannier,
{\it $S$-unit points on analytic hypersurfaces}, 
Ann. Sci. \'Ecole Norm. Sup. 38 (2005), 76--92.

\bibitem{CoZa12}
P. Corvaja and U. Zannier,
{\it Finiteness of odd perfect powers with four nonzero binary digits},
Ann. Inst. Fourier (Grenoble) 63 (2013), 715--731.

\bibitem{GrVi13}
S. S. Gross and A. F. Vincent,
{\it On the factorization of $f(n)$ for $f(x)$ in $\Z[x]$}, 
Int. J. Number Theory 9 (2013), 1225--1236.

\bibitem{GyYu06}
K. Gy\H ory and K. Yu, 
{\it Bounds for the solutions of $S$-unit equations and decomposable form equations}, 
Acta Arith.  123 (2006), 9--41. 

\bibitem{Matv00} 
E.\ M.\ Matveev,
{\it An explicit lower bound for a homogeneous rational linear form
in logarithms of algebraic numbers.\ II},
Izv.\ Ross.\ Acad.\ Nauk Ser.\ Mat.\  {64}  (2000),  125--180 (in Russian); 
English translation in Izv.\ Math.\  {64} (2000),  1217--1269.

\bibitem{MauRi10}
Ch. Mauduit and J. Rivat,
{\it Sur un probl\`eme de Gelfond : la somme des chiffres des
nombres premiers}, 
Ann. of Math. 171 (2010), 1591--1646.

\bibitem{May16}
J. Maynard,
{\it Primes with restricted digits}.
Preprint. Available at: 

 {\tt http://arxiv.org/pdf/1604.01041v1.pdf}

\bibitem{Schm70a}
{W. M. Schmidt},
{\it Simultaneous approximations to algebraic numbers by rationals},
Acta Math. 125 (1970), 189--201.

\bibitem{Schm72}
W. M. Schmidt,
{\it Norm form equations},
Ann. of Math. 96 (1972), 526--551.

\bibitem{SchmLN}
W. M. Schmidt,
Diophantine Approximation, Lecture Notes in Math. 
{785}, Springer, Berlin, 1980.

\bibitem{Ste08c}
C. L. Stewart,
{\it On the greatest square-free factor of terms of a linear recurrence sequence}. 
In: Diophantine equations, 257--264, 
Tata Inst. Fund. Res. Stud. Math., 20, Tata Inst. Fund. Res., Mumbai, 2008.

\bibitem{Ste13b}
C. L. Stewart,
{\it On prime factors of terms of linear recurrence sequences}. 
In: Number theory and related fields, 341--359, 
Springer Proc. Math. Stat., 43, Springer, New York, 2013.

\bibitem{Yu07}
K. Yu, 
{\it $p$-adic logarithmic forms and group varieties. III}, 
Forum Math. 19 (2007),187--280.

\vskip 1cm

\noi Yann Bugeaud

\noi Institut de Recherche Math\'ematique Avanc\'ee, U.M.R. 7501

\noi Universit\'e de Strasbourg et C.N.R.S.

\noi 7, rue Ren\'e Descartes

\noi 67084 STRASBOURG \ \ (France)

\vskip 2mm

\noi e-mail : {\tt bugeaud@math.unistra.fr}

\bye